\documentclass[12pt]{article}

\usepackage{amsmath}
\usepackage{amssymb}
\usepackage{eucal}

\usepackage[english]{babel}

\begin{document}

\baselineskip 16pt

\title{Finite soluble groups with nilpotent wide subgroups}


\author{V.\,S.~Monakhov and I.\,L.~Sokhor}


\maketitle

\begin{abstract}

A subgroup $H$ of a finite group $G$ is wide
if each prime divisor of the order of~$G$
divides the order of~$H$. We obtain the description
of finite soluble groups with no wide subgroups.
We also prove that a finite soluble group with
nilpotent wide subgroups has the quotient group
by its hypercenter with no wide subgroups.

\end{abstract}




{\small }

\section{Introduction}

\bigskip

All groups in this paper are finite.
All notations and terminology are standard.
The reader is referred to~\cite{Mon, Hup} if necessary.

Let $G$ be a group. We use $\pi(G)$ to denote the set
of all prime devisors of~$|G|$. By  $|\pi(G)|$ we denote
a number of different prime devisors of~$|G|$.
A subgroup $H$ of a group $G$ is said to be wide
if $\pi (H)=\pi (G)$.

Let $k$ be a positive integer.
A group $G$ is called $k$-primary if $|\pi(G)|=k$.
If $|\pi(G)|=1$ or $|\pi(G)|=2$, then $G$ is said to be
primary or biprimary, respectively.
A group $G$ is quasi-$k$-primary if $|\pi (G)|>k$ and
$|\pi (M)|\le k$ for every maximal subgroup~$M$ of~$G$.
Obviously, quasi-$k$-primary groups have
no wide subgroups.

A quasi-$1$\nobreakdash-\hspace{0pt}primary group is
also called quasiprimary, and
a quasi-$2$\nobreakdash-\hspace{0pt}primary group
is also called quasibiprimary~\cite{Liv}.
It is clear that the order of a nilpotent quasiprimary group
is equal to ~$pq$, where $p$ and $q$ are primes.
A nonnilpotent quasiprimary group~$G$
can be represented as the semidirect product $G=[E_{p^a}]Z_q$ of
a normal elementary abelian group~$E_{p^a}$, $|E_{p^a}|=p^a$,
by a cyclic group~$Z_q$, $|Z_q|=q$, where~$a$ is the exponent of~$p$.
This follows from Schmidt theorem~\cite{Shm} on groups all of whose
subgroups are nilpotent.

S.\,S.~Levischenko investigated quasibiprimary groups~\cite{Liv}.
A soluble quasibiprimary group~$G$ can be represented
as the semidirect product~$[P]M$  of an elementary abelian Sylow
$p$\nobreakdash-\hspace{0pt}subgroup~$P$ and a quasiprimary subgroup~$M$,
which is a maximal subgroup of ~$G$~\cite[Theorem~3.1]{Liv}.
In an insoluble quasibiprimary group~$G$ the Frattini subgroup $\Phi (G)$
is primary~\cite[Theorem 2.2]{Liv}, the quotient group~$G/\Phi (G)$
is a simple group, and all such simple groups are enumerated~\cite[Theorem~2.1]{Liv}.

It is natural to study the structure of quasi-$k$\nobreakdash-\hspace{0pt}primary
groups for any positive integer $k$. Most of simple groups is quasi-$k$\nobreakdash-\hspace{0pt}primary.
Simple groups that are not quasi-$k$\nobreakdash-\hspace{0pt}primary are enumerated
in~\cite[3.8]{ZW2005}. It is satisfied the question of V.\,S.~Monakhov
in Kourovka notebook~\cite[11.64]{Kour}.

We obtain the description of soluble quasi-$k$\nobreakdash-\hspace{0pt}primary groups
and soluble groups with nilpotent wide subgroups.


\section{Preliminaries}


{\bf Lemma~1.}
{\sl A soluble quasi-$k$\nobreakdash-\hspace{0pt}primary group is $(k+1)$\nobreakdash-\hspace{0pt}primary.}

\medskip

{\sc Proof.}
Let $G$ be a soluble quasi-$k$\nobreakdash-\hspace{0pt}primary
group and  $M$ be a maximal subgroup of $G$.
Then $|\pi(G)|\geq k+1$ and $|\pi(M)|\leq k$.
In a soluble group maximal subgroups have prime indices.
Therefore
$$
|G:M|=p^{\alpha}, \ p\in \pi(G), \ \alpha \in \mathbb{N}.
$$
Since $|G|=|M|\cdot |G:M|$, we have $|\pi(G)|\leq k+1$, and $|\pi(G)|=k+1$.
Lemma~1 is proved.

\medskip

{\bf Lemma~2.}
{\sl If $G$ is a soluble quasi-$k$\nobreakdash-\hspace{0pt}primary
group and $N$ is its normal Hall subgroup, then
$G/N$ is a quasi-$l$\nobreakdash-\hspace{0pt}primary group,
where $l=k-|\pi(N)|$.}

\medskip

{\bf Proof.}
By Lemma~1,  $G$ is $(k+1)$\nobreakdash-\hspace{0pt}primary,
hence
$$
|\pi(G/N)|=|\pi(G)|-|\pi(N)|=k+1-|\pi(N)|=l+1.
$$
Let $M/N$ be a maximal subgroup of~$G/N$.
Then $M$ is a maximal subgroup of~$G$.
Therefore $|\pi(M)|<|\pi(G)|$.
As $N$ is a Hall subgroup of~$G$, we obtain
$$
|\pi(M/N)|=|\pi(M)|-|\pi(N)|<|\pi(G)|-|\pi(N)|=|\pi (G/N)|=l+1.
$$
Thus $G/N$ is $(l+1)$\nobreakdash-\hspace{0pt}primary,
and any maximal subgroup of~$G$ is not more than
$l$\nobreakdash-\hspace{0pt}primary.
Consequently, $G/N$ is quasi-$l$\nobreakdash-\hspace{0pt}primary.
Lemma~2 is proved.

\medskip

{\bf Lemma~3.} (\cite[IV.2.11]{Hup}) {\it
If all Sylow subgroups of a group $G$ are cyclic,
then the derived subgroup $G^\prime$ is a cyclic Hall subgroup and
the quotient group $G/G^\prime$ is cyclic.}


\section{The structure of soluble quasi-$k$-primary groups}


{\bf Theorem~1.} {\sl Let $G$ be a soluble group.
Then the following statements are equivalent.}

(1)~{\sl $G$ is a quasi-$k$\nobreakdash-\hspace{0pt}primary group.}

(2)~{\sl Every normal subgroup of~$G$ is a Hall subgroup.}

(3)~{\sl Every maximal subgroup of~$G$ is a Hall subgroup.}

(4)~{\sl $G=[N]M$, where $N$ is a minimal normal and Sylow subgroup of~$G$,
$M$ is a  quasi-$(k-1)$\nobreakdash-\hspace{0pt}primary and maximal subgroup.}

\medskip

{\bf Proof.}
Check that $(1)$ implies $(2)$.
Let $G$ be a soluble quasi-$k$\nobreakdash-\hspace{0pt}primary group
and $N$ be a normal subgroup of~$G$, $\tau =\pi(G)\setminus \pi(N)$.
Suppose that $N$ is not a Hall subgroup of~$G$.
Since $G$ is a quasi-$k$\nobreakdash-\hspace{0pt}primary group
and $N$ is its proper subgroup, we have $|\pi(N)|\leq k$
and $\tau\neq \emptyset$.
As $G$ is a soluble group, in $G$ there is
$\tau$\nobreakdash-\hspace{0pt}Hall subgroup $M$.
Now, $(|M|,|N|)=1$, so $M\cap N=1$ and $NM=[N]M<G.$
But
$$
\pi([N]M)=\pi(N)\cup \pi(M)=\pi(N)\cup \tau=\pi(G), \ |\pi(G)|=k+1,
$$
that is, in a quasi-$k$\nobreakdash-\hspace{0pt}primary group $G$
there is a proper $(k+1)$\nobreakdash-\hspace{0pt}primary subgroup~$[N]M$,
a contradiction.  Thus $N$  is a Hall subgroup of $G$.

Now we prove that $(3)$ follows from $(2)$.
Assume that in a soluble group $G$
every normal subgroup is a Hall subgroup.
Let $M$ be a maximal subgroup and
$N$ a minimal normal subgroup  of $G$.
Then $N$ is a Sylow subgroup of $G$.
If $N$ does not belong to~$M$, then $G=NM$.
And since $N$ is abelian, $N\cap M=1$ and subgroup $M$ is a Hall subgroup.
Let $N\subseteq M$. Then $M/N$ is a maximal subgroup of~$G/N$.
It is clear that every normal subgroup of $G/N$ is a Hall subgroup.
Then by induction  $M/N$ is a Hall subgroup of $G/N$.
Therefore $M$ is also a Hall subgroup of~$G$.

Check that $(3)$ implies $(1)$.
Assume that every maximal subgroup of a soluble group~$G$
is a Hall subgroup. Let $M$ be a maximal subgroup of~$G$
and $|\pi(G)|=k+1$. Since $M$ is a Hall subgroup,
we have $|\pi(M)|=|\pi(G)|-1=k$, and
$G$ is a quasi-$k$\nobreakdash-\hspace{0pt}primary group.

Thus, $(1)$, $(2)$ and $(3)$ are equivalent.

Check that $(1)$ implies $(4)$.
Let $G$ be a soluble quasi-$k$\nobreakdash-\hspace{0pt}primary
group and $N$ its minimal normal subgroup.
In view of $(2)$, $N$ is a Sylow $p$\nobreakdash-\hspace{0pt}subgroup of~$G$.
By  Schur–Zassenhaus theorem~\cite[Theorem~4.32]{Mon},
in $G$ there is a subgroup $M$ such that $G=[N]M$.
Applying Lemma~1, $G$ is $(k+1)$\nobreakdash-\hspace{0pt}primary.
Hence $M$ is $k$\nobreakdash-\hspace{0pt}primary.
Suppose that in $M$ there is a proper
$k$\nobreakdash-\hspace{0pt}primary subgroup $M_1$.
Then in $G$ there is a proper
$(k+1)$\nobreakdash-\hspace{0pt}primary subgroup $[N]M_1$,
a contradiction. Consequently, $M$ is
quasi-$(k-1)$\nobreakdash-\hspace{0pt}primary.

Finely we prove that $(1)$ follows from $(4)$.
Let $G=[N]M$, where $N$ is a minimal normal and
Sylow subgroup of $G$, $M$ is its quasi-$(k-1)$\nobreakdash-\hspace{0pt}primary
and maximal subgroup. Then by Lemma~1,
$$
|\pi(M)|=k, \ |\pi(G)|=|\pi([N]M)|=|\pi(N)|+|\pi(M)|=k+1.
$$
Show that every maximal subgroup $K$ of $G$ is
not more than $k$\nobreakdash-\hspace{0pt}primary.
If $NK=G$, then $N\cap K=1$ and $K\simeq M$.
Since $M$ is a quasi-$(k-1)$\nobreakdash-\hspace{0pt}primary
soluble subgroup of $G$, using Lemma~1, we obtain that
$M$ is $k$\nobreakdash-\hspace{0pt}primary.
And so $K$ is also a $k$\nobreakdash-\hspace{0pt}primary group.
Assume that $NK$ is a proper subgroup of~$G$.
Then $N\subseteq K$ and $K=N(K\cap M)$ by Dedekind identity.
As $M$ is quasi-$(k-1)$\nobreakdash-\hspace{0pt}primary,
we obtain that $L\cap M$ is not more than
$(k-1)$\nobreakdash-\hspace{0pt}primary.
Hence $K=[N](K\cap M)$ is not more than
$k$\nobreakdash-\hspace{0pt}primary.
Thus $(4)$ and $(1)$ are equivalent.
Theorem~1 is proved.

\medskip

Note that $\pi$\nobreakdash-\hspace{0pt}soluble groups
with certain Hall maximal subgroups are  investigated in~\cite{Mon08}.

\medskip

{\bf Corollary 1.1.}
{\sl In a soluble quasi-$k$\nobreakdash-\hspace{0pt}primary group
the Frattini subgroup is trivial.}

\medskip

{\bf Proof.}
The Frattini subgroup of any group can never be
a Hall subgroup~\cite[Theorem~4.33]{Mon}.
It remains only to use Statement~$(2)$ of Theorem~1.

\medskip

{\bf Corollary 1.2.}
{\sl In a soluble quasi-$k$\nobreakdash-\hspace{0pt}primary group
the Fitting subgroup is a Hall subgroup  and every its Sylow subgroup
is a minimal normal subgroup.}

\medskip

{\bf Proof.}
By Theorem~1\,$(2)$ the Fitting subgroup of a soluble
quasi-$k$\nobreakdash-\hspace{0pt}primary group $G$ is a Hall subgroup.
And so every its Sylow subgroup $P$ is a Sylow subgroup of $G$.
Moreover, since the Fitting subgroup is nilpotent,
$P$ is a characteristic subgroup. Consequently,
$P$ is normal in $G$. At the same time by Theorem~1\,$(2)$,
a minimal normal subgroup of a soluble
quasi-$k$\nobreakdash-\hspace{0pt}primary group
is a Sylow subgroup. Corollary~1.2 is proved.

\medskip

{\bf Corollary 1.3.} {\sl Let $N$ be a normal subgroup of a group~$G$.
If  $G$ is a soluble quasi-$k$\nobreakdash-\hspace{0pt}primary group,
then the quotient group $G/N$ is a soluble
quasi-$l$\nobreakdash-\hspace{0pt}primary group,
where~$l=k-|\pi(N)|$.}

\medskip

{\bf Proof.} It follows from Theorem~1\,$(2)$ and Lemma~2.

\medskip

If a group $G$ has a normal series whose factors are isomorphic
to Sylow subgroups, then we say  $G$ has a Sylow tower.

\medskip

{\bf Corollary~1.4.}
{\sl A soluble quasi-$k$\nobreakdash-\hspace{0pt}primary group has a Sylow tower.}

\medskip

{\bf Proof.} We use induction on $k$.
Let $G$ be a soluble quasi-$k$\nobreakdash-\hspace{0pt}primary group.
Then by Statement~$(2)$ and $(4)$ of Theorem~1,
$G$ can be represented as $G=[P_1]M_1$,
where $P_1$ is a minimal normal and Sylow subgroup of $G$,
$M_1$ is a quasi-$(k-1)$\nobreakdash-\hspace{0pt}primary and
maximal subgroup. By the induction hypothesis,
$M_1$ has a Sylow tower. Hence $G$ also has a Sylow tower.
Corollary~1.4 is proved.

\medskip

A positive integer $n$ is said to be squarefree
if $p^2$ does not divide $n$ for all primes $p$.
A group is supersoluble if all its chief factors
are of prime orders.

\medskip

{\bf Corollary~1.5.}
{\sl The order of a group $G$ is squarefree if and only if
$G$ is a supersoluble quasi-$k$\nobreakdash-\hspace{0pt}primary
group, where~$k=|\pi (G)|-1$.
In particular, a supersoluble quasi-$k$\nobreakdash-\hspace{0pt}primary
group is metacyclic.}

\medskip

{\bf Proof.}
Let the order of a group $G$ be squarefree and $k=|\pi (G)|-1$.
Then all Sylow subgroups of~$G$ are cyclic  and
$G$ is supersoluble and metacyclic by Lemma~3.
It is clear that $|\pi (X)|<|\pi (G)|$
for every proper subgroup~$X$ of $G$,
that is, $G$ is quasi-$k$\nobreakdash-\hspace{0pt}primary.

Converse, let $G$ be a supersoluble
quasi-$k$\nobreakdash-\hspace{0pt}primary group.
Apply induction on~$k$. By Theorem~1\,$(4)$,
$G=[N]M$, where $N$ is a minimal normal and Sylow subgroup of $G$,
$M$ is a quasi-$(k-1)$\nobreakdash-\hspace{0pt}primary
and maximal subgroup. In view of~\cite[Theorem~4.48]{Mon},
a minimal normal subgroup of a supersoluble group
is of prime order, hence $|N|=p$, $p\in \pi (G)$.
A subgroup $M$ is supersoluble and
quasi-$(k-1)$\nobreakdash-\hspace{0pt}primary.
By  the induction hypotheses, the order of $M$
is squarefree. Hence the order of~$G$ is also squarefree.
And so~$G$ is metacyclic by Lemma~3.
Corollary~1.5 is proved.

\medskip

{\bf Corollary~1.6.}
{\sl The derived length of a soluble
quasi-$k$\nobreakdash-\hspace{0pt}primary group~$G$
does not exceed
$$
\min \{|\pi (G)|, \ \max \{1+a_i\mid i=1,2,\ldots ,t\} \},
$$
where $|F(G)|=p_1^{a_1}p_2^{a_2}\ldots p_t^{a_t}$.}

\medskip

{\bf Proof.}
Let $d(G)$ and $r(G)$ be the derived length and  the rank
of a soluble quasi-$k$\nobreakdash-\hspace{0pt}primary group~$G$,
respectively. By Theorem~1\,$(4)$,
all Sylow subgroups of $G$ are elementary abelian
and the length of chief series of $G$ equals $|\pi (G)|$.
Hence $d(G)\le |\pi (G)|$.
Soluble groups with abelian Sylow subgroups
are considered in paragraph~VI.14 \cite{Hup}.
Therefore all statements applying to
$A$\nobreakdash-\hspace{0pt}groups
(soluble groups with abelian Sylow subgroups)
of this paragraph
is correct for quasi-$k$\nobreakdash-\hspace{0pt}primary groups.
In particular, $d(G)\le 1+r(G)$
by~\cite[VI.14.18]{Hup}, and $r(G)\le u$ by~\cite[VI.14.31]{Hup},
where $u$ is the maximal number of
generating Sylow subgroups in the Fitting subgroup~$F(G)$.
Since Sylow subgroups of $G$ are elementary abelian,
we have $u=\max \{a_i\mid i=1,\ldots ,t \}$.
It follows that
$$
d(G)\le \min \{|\pi (G)|,\max \{1+a_i\mid i=1,\ldots ,t\}\}.
$$
Corollary~1.6 is proved.

\medskip

Substituting  $k=2$ in Theorem~1,
we obtain the result of S.\,S. Levischenko.

\medskip

{\bf Corollary~1.7.} (\cite[Theorem 3.1]{Liv})
{\sl A soluble group~$G$ is quasibiprimary if and only if
$G$ is equal to the semidirect product~$[P]M$ of its elementary abelian Sylow
$p$\nobreakdash-\hspace{0pt}subgroup~$P$ and quasiprimary subgroup~$M$,
which is also a maximal subgroup of~$G$. }


\section{Soluble groups with nilpotent wide subgroups}


Let $G$ be a nontrivial group,
$Z_0(G)=1$,
$Z_1(G)=Z(G)$, $Z_2(G)/Z_1(G)=Z(G/Z_1(G))$,
\ldots, $Z_i(G)/Z_{i-1}(G)=Z(G/Z_{i-1}(G))$, \ldots.
Then the subgroup
$Z_{\infty}(G)=\bigcup_{i=0}^{\infty} Z_i(G)$
is called the hypercenter of $G$.

Obviously, $Z(G/Z_{\infty}(G))=1$.

\medskip

{\bf Theorem~2.}
{\sl If in a soluble group $G$ every maximal subgroup~$M$ 
such that $\pi (M)=\pi (G)$ is nilpotent,
then the quotient group $G/Z_{\infty}(G)$ is
quasi-$k$\nobreakdash-\hspace{0pt}primary,
where $k=|\pi(G/Z_{\infty}(G))|-1$.}

\medskip

{\bf Proof.}
Suppose that the quotient group $\overline {G}=G/Z_{\infty}(G)$
is not quasi-$k$\nobreakdash-\hspace{0pt}primary, where
$k=|\pi(G/Z_{\infty}(G))|-1$. Then in $\overline {G}$
there is a maximal subgroup $\overline {M}=M/Z_{\infty}(G)$
such that $\pi(\overline {M})=\pi(\overline {G})$.
Hence $M$ is a maximal subgroup of $G$ and $\pi(M)=\pi(G)$.
By the hypotheses, $M$ is nilpotent.
Therefore $\overline {M}$ is also nilpotent.
Since $\overline {M}$ is a maximal subgroup
of a soluble group $\overline {G}$, we obtain that
$|\overline {G}:\overline {M}|=p^{\alpha}$
for some $p\in \pi(G)$ and positive integer~$\alpha$.
Let $\overline {M}_p$ be a Sylow
$p$\nobreakdash-\hspace{0pt}subgroup of~$\overline {M}$,
$\overline {G}_p$ be a Sylow $p$\nobreakdash-\hspace{0pt}subgroup
of~$\overline {G}$ containing $\overline {M}_p$.
Then $\overline {M}_p$ is a proper subgroup of
$\overline {G}_p$. Therefore $\overline {M}_p$
is normal in $\overline {G}$. Now there is
a nontrivial element $\overline {x}$
of $\overline {M}_p \cap Z(\overline {G}_p)$
such that it belongs to the center of~$\overline {G}$.
This contradicts the fact that $Z(G/Z_{\infty}(G))=1$.
Theorem~2 is proved.

\medskip

{\bf Corollary~2.1.}
{\sl If all wide subgroups of a group $G$ are nilpotent,
then $G/Z_{\infty}(G)$ is quasi-$k$\nobreakdash-\hspace{0pt}primary,
where $k=|\pi(G/Z_{\infty}(G))|-1$.}

\medskip

The converse of Theorem~2 is not true.
For example, let $G=S_3\times Z_6$.
Here $S_3$ is the symmetric group of order~3,
$Z_6$  is the cyclic group of order 6.
Then $Z_{\infty}(G)=Z_6$, and so $G/Z_{\infty}(G)\simeq S_3$
is quasiprimary, but the wide maximal subgroup $M=S_3\times Z_2$
is nonnilpotent.

\medskip

The following question becomes natural.

\emph{What is the structure of a finite soluble group
all of whose wide subgroups are supersoluble?}

\end{document}